\documentclass[12pt]{amsart}

\newtheorem{proposition}{Proposition}[section]
\newtheorem{lemma}[proposition]{Lemma}
\newtheorem{corollary}[proposition]{Corollary}

\theoremstyle{definition}

\theoremstyle{remark}
\newtheorem{remark}[proposition]{Remark}
\newtheorem{remarks}[proposition]{Remarks}

\newcommand{\selabel}[1]{\label{se:#1}}
\newcommand{\seref}[1]{Section~\ref{se:#1}}
\newcommand{\lelabel}[1]{\label{le:#1}}
\newcommand{\leref}[1]{Lemma~\ref{le:#1}}
\newcommand{\prlabel}[1]{\label{pr:#1}}
\newcommand{\prref}[1]{Proposition~\ref{pr:#1}}
\newcommand{\colabel}[1]{\label{co:#1}}
\newcommand{\coref}[1]{Corollary~\ref{co:#1}}
\newcommand{\relabel}[1]{\label{re:#1}}
\newcommand{\reref}[1]{Remark~\ref{re:#1}}

\newcommand{\eqlabel}[1]{\label{eq:#1}}
\newcommand{\equref}[1]{(\ref{eq:#1})}

\newcommand{\Hom}{{\rm Hom}}

\newcommand{\Ext}{{\rm Ext}}

\def\ot{\otimes}

\def\text#1{{\rm {\rm #1}}}

\begin{document}
\title[Semisimplicity of Yetter-Drinfeld modules]{Semisimplicity of the
categories of Yetter-Drinfeld modules and Long dimodules}
\author{S. Caenepeel}
\address{Faculty of Applied Sciences,
Vrije Universiteit Brussel, VUB, B-1050 Brussels, Belgium}
\email{scaenepe@vub.ac.be}
\urladdr{http://homepages.vub.ac.be/~scaenepe/}
\author{T. Gu\'ed\'enon}
\address{Faculty of Applied Sciences,
Vrije Universiteit Brussel, VUB, B-1050 Brussels, Belgium}
\email{guedenon@caramail.com}
\thanks{Research supported by the project G.0278.01 ``Construction
and applications of non-commutative geometry: from algebra to physics"
from FWO Vlaanderen}
\begin{abstract}
Let $k$ be a field, and $H$ a Hopf algebra with bijective antipode. 
If $H$ is commutative, noetherian, semisimple and cosemisimple, then
the category ${}_{H}{\mathcal {YD}}^H$ of Yetter-Drinfeld modules is semisimple.
We also prove a similar statement for the category of Long dimodules,
without the assumption that $H$ is commutative.
\end{abstract}
\maketitle

\section*{Introduction}\selabel{0}
Let $H$ be a Hopf algebra at the same time acting and coacting on 
a vector space $M$. We can impose various compatibility relations
between the action and coaction, leading to different notions
of Hopf modules. Hopf modules are already considered by Sweedler
\cite{Sweedler}, and they have to satisfy the relation
$$\rho(hm)=\Delta(h)\rho(m)=h_1m_0\ot h_2m_1$$
One can also require that the $H$-coaction is $H$-linear:
$$\rho(hm)=h\rho(m)=hm_0\ot m_1$$
A module satisfying this condition is called a Long dimodule.
Long dimodules are the  building stones of the Brauer-Long group,
in the case where the Hopf algebra $H$ is commutative, cocommutative
and faithfully projective (see \cite{Long}, and \cite{Caenepeel98}
for a detailed discussion). Long dimodules are also connected to a
non-linear equation (see \cite{Militaru97a}).\\
Another - at first sight complicated and artificial - compatibility
relation is the following:
$$h_1m_0 \otimes h_2m_1=(h_2m)_0 \otimes (h_2m)_1h_1$$
A module that satisfies it is called a Yetter-Drinfeld module. 
There is a close connection between Yetter-Drinfeld modules and the Drinfeld
double (see \cite{Drinfeld}): if $H$ is finitely generated projective,
then the category of Yetter-Drinfeld modules is isomorphic to the
category of modules over the Drinfeld double. Yetter-Drinfeld modules have
been studied intensively by several authors over the passed fifteen years,
see for example \cite{CaenepeelMZ02}, \cite{Lambe}, \cite{Majid},
\cite{Radford}, this list is far from exhaustive. One of the important features is the
fact that the category of Yetter-Drinfeld modules is braided monoidal. As Long
dimodules, Yetter-Drinfeld modules are related to a non-linear equation, the
quantum Yang-Baxter equation (see e.g. \cite{Kassel}). If $H$ is
commutative and cocommutative, then Yetter-Drinfeld modules coincide with
Long dimodules.\\
In this note, we give sufficient conditions for the categories of Yetter-Drinfeld
modules and Long dimodules to be semisimple (\seref{3}) and we study
projective and injective dimension in these categories. Our main result
is that the category of Yetter-Drinfeld modules is semisimple if
$H$ is a commutative, noetherian, semisimple and cosemisimple Hopf algebra
over a field $k$. The same is true for the category of Long dimodules,
without the assumption that $H$ is commutative.\\
For generalities on Hopf algebras, we refer the reader to \cite{DascalescuNR},
\cite{Montgomery}, \cite{Sweedler}. For a detailed study of Hopf modules
and their generalizations, we refer to \cite{CaenepeelMZ02}.

\section{Preliminary Results}\selabel{1}
Let $k$ be a commutative ring, and $H$ a faithfully flat Hopf algebra
with bijective antipode $S$. Unadorned $\ot$ and $\Hom$ will be over $k$.
We will use the Sweedler-Heyneman notation for comultiplication and
coaction: for $h\in H$, we write
$$\Delta(h)= h_1\ot h_2$$
(summation implicitly understood), and for a right $H$-comodule $(M,\rho_M)$
and $m\in M$, we write
$$\rho_M(m)= m_0\ot m_1$$
${}_{H}{\mathcal M}$ and ${\mathcal M}^H$ will be the categories of
respectively left $H$-modules and left $H$-linear maps, and right
$H$-comodules and right $H$-colinear maps. If $M$ and $N$ are right
$H$-comodules, then we denote the $k$-module consisting of right $H$-colinear
maps from $M$ to $N$ by $\Hom^H(M,N)$.
$$M^{{\rm co}H}=\{ m \in M~|~{\rho}_M (m) = m \otimes 1\}$$
is called the $k$-submodule of coinvariants of $M$. Observe that
$H^{{\rm co}H}=k$.\\
Suppose that a $k$-vector space $M$ is at the same time a left $H$-module
and a right $H$-comodule. Recall that $M$ is called a left-right
Yetter-Drinfeld module if
$$h_1m_0 \otimes h_2m_1=(h_2m)_0 \otimes (h_2m)_1h_1$$
or, equivalently,
$$\rho(hm)=h_2m_0 \otimes h_3m_1S^{-1}(h_1)$$ 
for all $m\in M$ and $h\in H$. $M$ is called a left-right Long dimodule if
$$\rho(hm)=hm_0\ot m_1$$
for all $m\in M$ and $h\in H$. If $H$ is commutative and cocommutative,
then a Long dimodule is the same as a Yetter-Drinfeld module. 
${}_{H}{\mathcal {YD}}^H$ and ${}_{H}{\mathcal L}^H$ will be the categories
of respectively Yetter-Drinfeld modules and Long dimodules, and $H$-linear
$H$-colinear maps. The $k$-module consisting of all $H$-linear $H$-colinear
maps between two Yetter-Drinfeld modules or two Long dimodules $M$ and $N$
will be denoted by ${}_H\Hom^H(M,N)$. If $H$ is finitely generated and
projective, then the category ${}_{H}{\mathcal {YD}}^H$ is isomorphic
to the category ${}_{D(H)}{\mathcal M}$, where $D(H)$ is the Drinfeld double
of $H$, and ${}_{H}{\mathcal L}^H$ is isomorphic to ${}_{H\ot H^*}{\mathcal M}$.\\
The functors
$$(-)^{{\rm co}H}:\ {}_{H}{\mathcal {YD}}^H\to {\mathcal M}~~{\rm and}~~
(-)^{{\rm co}H}:\ {}_{H}{\mathcal {L}}^H\to {\mathcal M}$$
are exact if
$$(-)^{{\rm co}H}:\ {\mathcal M}^H\to {\mathcal M}$$
is exact. This is the case if $H$ is cosemisimple and $k$ is a field.

\begin{lemma}\lelabel{1.1}
\begin{enumerate}
\item
Let $M$ and $N$ be objects of ${}_{H}{\mathcal {YD}}^H$. Then  $M \otimes N$ is an
object of $_{H}{\mathcal {YD}}^H$; the $H$-action and $H$-coaction are given by 
$$h(m \otimes n)=h_1m \otimes h_2n \quad \hbox{and} \quad \rho(m \otimes n)=m_0 \otimes n_0
\otimes n_1m_1$$ 
\item
Let $M$ and $N$ be objects of ${}_{H}{\mathcal L}^H$. Then  $M \otimes N$ is an
object of $_{H}{\mathcal L}^H$; the $H$-action and $H$-coaction are given by
$$h(m \otimes n)=h_1m \otimes h_2n \quad \hbox{and} \quad \rho(m \otimes n)=m_0 \otimes n_0  
\otimes m_1n_1$$ 
\item
For any $H$-comodule $N$, $H \otimes N$ is an object of $_{H}{\mathcal {YD}}^H$ via the
following structures
$$h(h' \otimes n)=hh' \otimes n \quad \hbox{and} \quad \rho(h \otimes n)=h_2 \otimes n_0
\otimes h_3n_1S^{-1}(h_1)$$
\item
For any $H$-comodule $N$, $H \otimes N$ is an object of $_{H}{\mathcal L}^H$ via the
following structures
$$h(h' \otimes n)=hh' \otimes n \quad \hbox{and} \quad \rho(h \otimes n)=h \otimes n_0 \otimes
n_1$$
\end{enumerate}
\end{lemma}

\begin{proof} This result is well-known, and the proof is a straightforward
computation. It may be found in \cite[p. 440]{Caenepeel98},
\cite[Prop. 12.1.2]{Caenepeel98}, \cite[Prop. 123]{CaenepeelMZ02},
and \cite[Sec. 7.2]{CaenepeelMZ02}.
\end{proof}

\begin{lemma}\lelabel{1.2}
\begin{enumerate} 
\item
Let $M$ and $N$ be in ${}_{H}{\mathcal {YD}}^H$. If $H$ is commutative, then
$M\otimes_HN$ is an object of ${}_{H}{\mathcal {YD}}^H$. The $H$-action and $H$-coaction
are given by 
$$h(m \otimes n)=hm \otimes n=m\otimes hn$$
and 
$${\rho}_{M\otimes_HN}(m\otimes n)=m_0 \otimes n_0\otimes n_1m_1$$ 
\item
Let $H$ be commutative. Let $M$ and $N$ be in ${}_{H}{\mathcal {YD}}^H$ with $M$
finitely generated projective in $_{H}{\mathcal M}$. Then 
\begin{enumerate}
\item
${}_H\Hom(M, N)\in {\mathcal M}^H$ and 
$${}_H\Hom^H(M, N)={}_H\Hom(M,N)^{coH}$$   
The coaction is defined by
$$\rho(f)=f_0 \otimes f_1 \in {}_H\Hom(M, N) \otimes H$$ 
if and only if 
\begin{equation}\eqlabel{1.2.1}
f_0(m) \otimes f_1=f(m_0)_0\otimes f(m_0)_1S(m_1)
\end{equation}
for all $m \in M$.
\item
${}_H\Hom(M, N)\in{}_{H}{\mathcal {YD}}^H$; the $H$-action is defined
by
$(hf)(m)=hf(m)=f(hm)$.
\end{enumerate}
\end{enumerate}
\end{lemma}

\begin{proof} 
1) It is clear that $M\otimes_HN$ is an $H$-module. An easy verification
shows that the $H$-coaction is well-defined on the tensor product over
$H$ and that the necessary associativity and counit properties are
satisfied, so that $M\otimes_HN$ is also an $H$-comodule.
$M\otimes_HN$ is a Yetter-Drinfeld module, since we have for
every $h\in H$ that
\begin{eqnarray*}
&&\hspace*{-2cm}
{\rho}_{M\otimes_H N}(hm\otimes n)=(hm)_0 \otimes n_0\otimes n_1(hm)_1\\
&=&h_2m_0 \otimes n_0 \otimes n_1h_3m_1S^{-1}(h_1)\\
&=&
h_2(m_0 \otimes n_0) \otimes h_3n_1m_1S^{-1}(h_1)\\
&=&h_2(m \otimes n)_0 \otimes h_3(m \otimes n)_1S^{-1}(h_1)
\end{eqnarray*}

2a) Let us define a map 
$$\pi :\ \Hom(M , N) \rightarrow \Hom(M , N \otimes H)$$
by 
$$\pi(f)(m)=f(m_0)_0 \otimes f(m_0)_1S(m_1)$$
Let $f$ be $H$-linear. Using the commutativity of $H$, we obtain
\begin{eqnarray*}
&&\hspace*{-2cm}
\pi(f)(hm)=f((hm)_0)_0 \otimes f((hm)_0)_1S((hm)_1)\\
&=&(h_2f(m_0))_0 \otimes (h_2f(m_0))_1S(h_3m_1S^{-1}(h_1))\\
&=&h_{3}f(m_0)_0 \otimes h_{4}f(m_0)_1S^{-1}(h_{2})h_1S(m_1)S(h_5)\\
&=&hf(m_0)_0 \otimes f(m_0)_1S(m_1)=h\pi(f)(m)
\end{eqnarray*}
so $\pi(f)$ is $H$-linear, and $\pi$ restricts to a map
$$\pi :\ {}_H\Hom(M , N) \rightarrow {}_H\Hom(M , N \otimes H)$$
Now $H$ is finitely generated and projective as an $H$-module, so we have
a natural isomorphism ${}_H\Hom(M , N \otimes H)\cong{}_H\Hom(M , N) \otimes H$,
and we obtain a map
$$\pi :\ {}_H\Hom(M , N) \rightarrow {}_H\Hom(M , N) \otimes H$$
with $\pi(f)=f_0\ot f_1$ if and only if
$$({\pi}(f))(m)=f_0(m) \otimes f_1=f(m_0)_0 \otimes f(m_0)_1S(m_1)$$
It is straightforward to show that $\pi$ makes ${}_H\Hom(M , N)$
a right $H$-comodule.
Now take
$f \in {}_H\Hom^H(M, N)$ and $m \in M$. Then
\begin{eqnarray*}
&&\hspace*{-2cm}
{\pi}(f)(m)=f_0(m) \otimes f_1 = f(m_0)_0 \otimes f(m_0)_1S(m_1)\\
&=&
f(m_{0}) \otimes m_{1}S(m_2)=f(m) \otimes 1= (f \otimes 1)(m)
\end{eqnarray*}
so $f$ is coinvariant. Conversely, take
$f \in {}_H\Hom(M, N)^{{\rm co}H}$. Then for every $m\in M$
$$f(m_0)_0\otimes f(m_0)_1S(m_1)=f_0(m) \otimes f_1 =f(m) \otimes 1$$ 
and 
$$f(m_{0})_{0} \otimes f(m_{0})_{1}S(m_1)m_2=f(m_0) \otimes m_1$$ 
and it follows that 
$${\rho}_N(f(m))={\rho}_N(f(m_0))\varepsilon (m_1)=f(m_0) \otimes m_1$$ 
and $f$ is $H$-colinear.\\
2b) Clearly ${}_H\Hom(M, N)$ is an $H$-module and, by a), it is an $H$-comodule. 
On the other hand, we have
\begin{eqnarray*}
&&\hspace*{-2cm}
((hf)_0 \otimes (hf)_1)(m)=((hf)(m_0))_0 \otimes ((hf)(m_0))_1S(m_1)\\
&=&
(h(f(m_0))_0 \otimes (h(f(m_0))_1S(m_1)\\
&=&h_2(f(m_0)_0) \otimes
h_3(f(m_0)_1)S^{-1}(h_1)S(m_1)\\
&=&h_2(f(m_0)_0) \otimes h_3(f(m_0)_1)S(m_1)s^{-1}(h_1)\\
&=&h_2(f_0(m)) \otimes h_3f_1S^{-1}(h_1)\\
&=&(h_2f_0 \otimes h_3f_1S^{-1}(h_1))(m)
\end{eqnarray*}
so  
${}_H\Hom(M, N)\in {}_{H}{\mathcal {YD}}^H$.
\end{proof}

\begin{remark}\relabel{1.3}
The results in \leref{1.2} remain true after we replace ${}_{H}{\mathcal YD}^H$
by ${}_{H}{\mathcal L}^H$. The $H$-coaction on
$M\otimes_H N$ is given by 
$${\rho}_{M\otimes_HN}(m\otimes n)=m_0 \otimes n_0\otimes m_1n_1$$ 
The $H$-coaction on ${}_H\Hom(M , N)$ is also defined by \equref{1.2.1}.
Part 2a) of \leref{1.2} then also holds if $H$
is noncommutative.
\end{remark} 

\begin{lemma}\lelabel{1.4}
Let $V$ be a $k$-module and $N$ an $H$-module.
\begin{enumerate}
\item ${}_H\Hom(H \otimes V , N)$ and $\Hom(V , N)$ are isomorphic
as $k$-modules.
\item If $V$ is projective as $k$-module, then $H \otimes V$ is projective in 
${}_H{\mathcal M}$.
\end{enumerate}
\end{lemma}

\begin{proof} 1) is well-known: the $k$-isomorphism 
$${}_H\Hom(H \otimes V, N) \rightarrow
\Hom(V, N)$$ is defined by $\phi (f)(v)=f(1 \otimes v)$.\\
2) follows immediately from (1).
\end{proof}

Let $V$ be an $H$-comodule which is finitely generated and projective as 
a $k$-module. By Lemmas \ref{le:1.1} and \ref{le:1.4},
$H \otimes V$ is an object in ${}_{H}{\mathcal {YD}}^H$ and in
${}_{H}{\mathcal L}^H$, and is finitely generated projective as an $H$-module. 
So if $N$
is an object of ${}_{H}{\mathcal {YD}}^H$ and if $H$ is commutative, then, by
\leref{1.2},
${}_H\Hom(H\otimes V, N)$ is an object in ${}_{H}{\mathcal {YD}}^H$. If $N$ is an object
of
${}_{H}{\mathcal L}^H$, then by \reref{1.3}, ${}_H\Hom(H\otimes V, N)$ is an object of
${\mathcal M}^H$; if furthermore $H$ is commutative, then ${}_H\Hom(H\otimes V, N)$ is
an object of ${}_{H}{\mathcal L}^H$.

\begin{lemma}\lelabel{1.5}
Let $H$ be commutative and $N\in {}_{H}{\mathcal {YD}}^H$. 
\begin{enumerate}
\item If $V$ is an $H$-comodule which is finitely generated and projective as
a $k$-module, then the
$H$-comodules ${}_H\Hom(H \otimes V, N)$ and $\Hom(V, N)$ are isomorphic.
\item  Let $k$ be a field and $V$ a finite-dimensional $H$-comodule that is projective
as an
$H$-comodule. Then $H\otimes V$ is a projective object of
 ${}_{H}{\mathcal{YD}}^H$.
\end{enumerate}
\end{lemma}

\begin{proof} 1) Consider the canonical $k$-isomorphism 
$$\phi :\ {}_H\Hom(H \otimes V, N) \rightarrow \Hom(V, N),~~
\phi(f)(v)=f(1 \otimes v)$$ 
$\phi$ is $H$-colinear since 
\begin{eqnarray*}
&&\hspace*{-2cm}
\phi(f)_0(v) \otimes \phi(f)_1=(\phi(f)(v_0))_0 \otimes (\phi(f)(v_0))_1S(v_1)\\
&=&
f(1\otimes v_0)_0 \otimes f(1 \otimes v_0)_1S(v_1)\\
&=&f({(1\otimes v)}_0)_0 \otimes f({(1\otimes v)}_0)_1S({(1 \otimes v)}_1)\\
&=&f_0(1 \otimes v) \otimes f_1\\
&=&(\phi(f_0))(v) \otimes f_1
\end{eqnarray*}
2) By 1) and \leref{1.2}, we have
\begin{eqnarray*}
&&\hspace*{-2cm}
{}_H\Hom^H(H \otimes V , N)\cong {}_H\Hom(H \otimes V , N)^{{\rm co}H}\\
&\cong &\Hom(V , N)^{{\rm co}H}\cong \Hom^H(V , N)
\end{eqnarray*}
\end{proof}

\leref{1.5} also holds with ${}_{H}{\mathcal{YD}}^H$ replaced by
${}_{H}{\mathcal{L}}^H$, and without the assumption that $H$ is commutative.

\begin{proposition}\prlabel{1.6} 
Let $k$ be a field. An object $M$ of $_{H}{\mathcal {YD}}^H$ or $_{H}{\mathcal
L}^H$ is finitely generated as an $H$-module if and only if there exists a
finite dimensional
$H$-comodule $V$ and an $H$-linear $H$-colinear epimorphism
$\pi:\ H\otimes V\longrightarrow M$.
\end{proposition} 

\begin{proof} If there exist a finite dimensional $H$-comodule $V$ and an 
epimorphism of
$H$-modules $\pi:\ H\otimes V\longrightarrow M$, then $H \otimes V$ is finitely
generated as an $H$-module and $M$ is a quotient of $H \otimes V$ in ${}_{H}{\mathcal
M}$, so $M$ is finitely generated in $_{H}{\mathcal M}$.\\
Suppose that $M$ is finitely generated as an $H$-module,
with generators $\{m_1,\cdots,m_n\}$. By \cite[5.1.1]{DascalescuNR}, there
exists a finite dimensional
$H$-subcomodule $V$ of $M$ containing $\{m_1,\cdots,m_n\}$ and the $k$-linear map 
$$\pi:\ H\otimes V \rightarrow M,~~\pi(h\otimes v)=hv$$
is an $H$-linear $H$-colinear epimorphism. 
\end{proof} 

Let $H^*$ be the linear dual of $H$. If $M$ and $N$ are $H$-comodules, then
 $\Hom_k(M ,N)$ is a left
$H^*$-module, with $H^*$-action
$$(h^*f)(m)=h^*(f(m_0)_1S(m_1))f(m_0)_0$$
(adapt the proof of \cite[Proposition 1.1]{Stefan}).

\begin{lemma}\lelabel{1.7}
Let $H$ be commutative. For $M,N\in {}_{H}{\mathcal {YD}}^H$, ${}_H\Hom(M ,
N)$ is a left $H^*$-submodule of $\Hom_k(M , N)$. 
\end{lemma}

\begin{proof} For all $\alpha \in H^*$, $f\in \Hom_H(M , N)$, $h\in H$ and $m\in M$,
we have
\begin{eqnarray*}
&&\hspace*{-2cm}
(\alpha f)(hm)=\alpha\Bigl(f((hm)_0)_1S((hm)_1)\Bigr)f((hm)_0)_0\\
&=&\alpha\Bigl(f(h_2m_0)_1S(h_3m_1S^{-1}(h_1))\Bigr)f(h_2m_0)_0\\
&=&\alpha\Bigl((h_2(f(m_0)))_1h_1S(m_1)S(h_3)\Bigl)(h_2f(m_0))_0\\
&=& \alpha\Bigl(h_{4}f(m_0)_1S^{-1}(h_{2})h_1S(m_1)S(h_5)\Bigr)h_{3}f(m_0)_0\\
&=&\alpha\bigl(f(m_0)_1S(m_1)\bigr)hf(m_0)_0\\
&=&h((\alpha f)(m))
\end{eqnarray*}
and it follows that $\alpha f$ is $H$-linear. Observe that we used
the commutativity of $H$.
\end{proof}

Recall that a left $H^*$-module $M$ is called rational if there exists a right
$H$-coaction on $M$ inducing the left $H^*$-action.

\begin{proposition}\prlabel{1.8}
Let $H$ be a commutative Hopf algebra over a field $k$. If $M,N\in
{}_{H}{\mathcal
{YD}}^H$ with $M$ finitely generated as $H$-module, then 
${}_H\Hom(M, N)\in {}_{H}{\mathcal {YD}}^H$.
\end{proposition}

\begin{proof} By \prref{1.6}, there exist a finite dimensional $H$-subcomodule $V$ of
$M$ and an $H$-linear $H$-colinear epimorphism $\pi:\ H\otimes V\longrightarrow
M$. So we obtain an injective $k$-linear map
$${}_H\Hom(\pi, N):\ {}_H\Hom(M, N)\rightarrow {}_H\Hom(H\otimes V, N)$$ 
For all $\alpha \in H^*$, $f\in {}_H\Hom(M , N)$, $h \in H$ and $v\in V$, we have
$\pi(h
\otimes v)=hv$,
$(1\otimes v)_0
\otimes (1\otimes v)_1=1\otimes v_0 \otimes v_1$ and 
\begin{eqnarray*}
&&\hspace*{-2cm}
(\alpha f)\circ \pi(1\otimes v)=(\alpha f)(v)=\alpha(f(v_0)_1S(v_1))f(v_0)_0\\
&=&\alpha(f(\pi(1\otimes v_0))_1S(v_1))f(\pi(1\otimes v_0))_0\\
&=&\alpha(f(\pi(1\otimes
v)_0))_1S((1\otimes v)_1)f(\pi((1\otimes v)_0))_0\\
&=&(\alpha(f\circ \pi))(1\otimes v)
\end{eqnarray*} 
This relation and the fact that $(\alpha f)\circ \pi$ and $\alpha(f\circ \pi)$ are $H$-linear imply that 
$((\alpha f)\circ \pi)(h\otimes v)=(\alpha(f\circ
\pi))(h\otimes v)$, and it follows that the map ${}_H\Hom(\pi, N)$ is
$H^*$-linear. By \leref{1.2},
${}_H\Hom(H\otimes V, N)$ is an $H$-comodule, and therefore a rational $H^*$-module. It
follows that ${}_H\Hom(M , N)$ is a rational $H^*$-module, being an
$H^*$-submodule of the rational $H^*$-module
${}_H\Hom(H\otimes V, N)$. This shows that ${}_H\Hom(M , N)$ is an $H$-comodule. By
\leref{1.2},  ${}_H\Hom(M , N)\in {}_{H}{{\mathcal {YD}}^H}$.
\end{proof}

\begin{remark} \relabel{1.9}
1) \leref{1.7} is still true if we replace ${}_{H}{{\mathcal {YD}}^H}$
by ${}_{H}{{\mathcal {L}}^H}$, without the assumption that $H$ is commutative.\\
2) We have the following Long dimodule version of \prref{1.8}: for
a (not necessarily commutative) Hopf algebra over a field $k$, and
$M,N\in {}_{H}{{\mathcal {L}}^H}$, with $M$ finitely generated as an $H$-module,
${}_H\Hom(M,N)$ $\in {\mathcal M}^H$.
\end{remark}

\section{Projective and injective dimension in the category 
of Yetter-Drinfeld modules}\lelabel{2}
 
\begin{lemma}\lelabel{2.1}
Let $H$ be commutative, and $M,N,P\in
{}_{H}{\mathcal {YD}}^H$, with $N$ finitely generated projective as an $H$-module. 
\begin{enumerate}
\item We have a $k$-isomorphism 
$${}_H\Hom^H( M , {}_H\Hom( N, P) )\cong {}_H\Hom^H(M\otimes _{H}N , P)$$ 
\item The functor 
$${}_H\Hom(N, -):\ {}_{H}{\mathcal {YD}}^H\to {}_{H}{\mathcal {YD}}^H$$ 
preserves injective objects.
\end{enumerate}
\end{lemma}

\begin{proof} 1) We have a natural isomorphism 
$$\phi :\ {}_H\Hom( M , {}_H\Hom( N, P))
\rightarrow {}_H\Hom( M\otimes {}_{H}N , P)$$
given by $\phi(f)(m\otimes n)=f(m)(n)$.
We will show that $\phi$ restricts to an isomorphism between
${}_H\Hom^H( M , {}_H\Hom( N, P) )$ and ${}_H\Hom^H(M\otimes _{H}N , P)$.
Take $f\in {}_H\Hom( M , {}_H\Hom( N, P))$ and $\phi(f)=g$. Then
$f$ is $H$-colinear if and only if
$$f(m_0)\ot m_1= f(m)_0\ot f(m)_1$$
for all $m\in M$.
Using \equref{1.2.1}, we find that this is equivalent to
$$f(m_0)(n)\ot m_1= f(m)_0(n)\ot f(m)_1=
f(m)(n_0)_0\ot f(m)(n_0)_1S(n_1)$$
for all $m\in M$ and $n\in N$, or
$$g(m_0\ot n)\ot m_1=g(m\ot n_0)_0\ot g(m\ot n_0)_1S(n_1)$$
which is equivalent to
$$g(m_0\ot n_0)\ot m_1n_1=g(m\ot n_0)_0\ot g(m\ot n)_1$$
and this equation means that $g$ is $H$-colinear.\\
2) If $I$ is an injective object of ${}_{H}{\mathcal {YD}}^H$,
then the functor 
$${}_H\Hom^H( -, I):\ {}_{H}{\mathcal {YD}}^H\to {}_{k}{\mathcal M}$$ 
is exact. On the other hand, $N$ is $H$-projective, hence the functor
$$(-)\otimes_{H}N:\ {}_{H}{\mathcal {YD}}^H\to {}_{H}{\mathcal {YD}}^H$$ 
is exact, and it follows from (1) that
$${}_H\Hom^H( - , {}_H\Hom( N, I)):\ {}_{H}{\mathcal {YD}}^H\to {}_{k}{\mathcal M}$$ 
is exact.
\end{proof}

If $k$ is a field, then the category of Yetter-Drinfeld modules 
${}_{H}{\mathcal {YD}}^H$ is Grothendieck, and every object has an injective
resolution. For every Yetter-Drinfeld module $M$, we can define the right derived
functors ${{}_H\Ext^H}^{i}( M , -)$ of
the covariant left exact functor 
$${}_H\Hom^H(M, -):\ {}_{H}{\mathcal {YD}}^H \rightarrow {}_k{\mathcal M}$$ 

\begin{proposition}\prlabel{2.2}
Let $H$ be a commutative Hopf algebra over a field $k$, and 
$M,N,P\in {}_{H}{\mathcal {YD}}^H$ with
$N$ finitely generated projective as an $H$-module. Then  
$${{}_H\Ext^H}^{i}( M ,  {}_H\Hom( N, P) ) \cong 
{{}_H\Ext^H}^i( M\otimes _{H}N , P)$$
\end{proposition}

\begin{proof} By the first part of \leref{2.1}, the functors 
$${}_H\Hom^H( M , {}_H\Hom( N, -))  \quad \hbox{and} 
\quad {}_H\Hom^H( M\otimes _{H}N ,-)$$  
coincide on ${}_{H}{\mathcal {YD}}^H$. By the projectivity and the finiteness
assumptions on $N$, the ${}_H\Hom(N, -)$ is an exact endofunctor of  
${}_{H}{\mathcal {YD}}^H$. By the second part of \leref{2.1}, it preserves the
injective objects of
${}_{H}{\mathcal {YD}}^H$. Thus the functor ${}_H\Hom(N, -)$ preserves injective
resolutions in ${}_{H}{\mathcal {YD}}^H$.
\end{proof}

In the following corollary, ${}_H{\rm pdim}^H(-)$ and ${}_H{\rm injdim}^H(-)$ denote
respectively the projective and injective dimension in the category ${}_{H}{{\mathcal
{YD}}^H}$.

\begin{corollary}\colabel{2.3}
Let $H$ be a commutative Hopf algebra over a field $k$, and 
$M,N,P\in {}_{H}{\mathcal {YD}}^H$ with
$N$ finitely generated projective as an $H$-module. Then 
\begin{enumerate}
\item ${}_H{\rm pdim}^H(M\otimes_HN)\leq {}_H{\rm pdim}^H(M)$.
\item ${}_H{\rm injdim}^H({}_H\Hom(N, P))\leq {}_H{\rm injdim}^H(P)$.
\end{enumerate}
\end{corollary}

\begin{remarks}\relabel{2.4}
1) Let $H$ be semisimple. Then the projectivity assumption in
\leref{2.1}, \prref{2.2} and \coref{2.3} is no longer needed.\\
2) If $k$ is a field, then ${}_{H}{\mathcal L}^H$ is a Grothendieck category 
with enough
injective objects, and every Long dimodule has an injective resolution. 
For every $M\in {}_{H}{\mathcal L}^H$, we can
then define the right derived functors ${{}_H\Ext^H}^{i}( M , -)$ of the covariant left
exact functor $${}_H\Hom^H(M, -):\ {}_{H}{\mathcal L}^H \rightarrow {}_k{\mathcal M}$$ 
All the results of this Section remain valid for ${}_{H}{\mathcal L}^H$. 
If $H$ is semisimple,
then the projectivity assumptions are not needed. 
\end{remarks} 

\section{Semisimplicity of the category of Yetter-Drinfeld modules}\selabel{3}
Throughout this Section, $k$ will be a field, and $H$ a commutative
Hopf algebra. Recall
that $M\in {}_{H}{\mathcal {YD}}^H$ is called simple if it has no proper
subobjects; a direct sum of simples is called semisimple. 
If every $M\in {}_{H}{\mathcal {YD}}^H$ is semisimple, then we call the
category ${}_{H}{\mathcal {YD}}^H$ semisimple.
We say that ${}_{H}{\mathcal {YD}}^H$ satisfies condition ($\dagger$)
if the following holds:\\
if $M\in {}_{H}{\mathcal {YD}}^H$ is finitely generated as a left
$H$-module, then ${}_H\Hom(M,-):\ {}_{H}{\mathcal {YD}}^H\to
{}_{H}{\mathcal {YD}}^H$ is exact.\\
By \prref{1.8}, ${}_H\Hom(M,N)\in {}_{H}{\mathcal {YD}}^H$ if
$H$ is commutative and $M$ is finitely generated as an $H$-module.
Also observe that 
${}_{H}{\mathcal {YD}}^H$ satisfies condition ($\dagger$) if $H$ is semisimple.

\begin{proposition}\prlabel{3.1}
Let $H$ be commutative. Assume that ${}_{H}{\mathcal {YD}}^H$ satisfies condition
($\dagger$) and that the functor 
$$(-)^{{\rm co}H}:\ {}_{H}{\mathcal {YD}}^H\to {\mathcal M}$$
 is exact. If $M\in {}_{H}{\mathcal {YD}}^H$ is finitely
generated as an $H$-module, then $M$ is a projective object in ${}_{H}{\mathcal
{YD}}^H$.
\end{proposition}

\begin{proof} 
We know that
$${}_H\Hom^H(M, -)\cong {}_H\Hom(M, -)^{{\rm co}H}$$
so ${}_H\Hom^H(M, -)$ is exact since it is isomorphic to the composition of two exact
functors.
\end{proof}

\begin{corollary}\colabel{3.2}
With the same assumptions as in \prref{3.1}, and with $H$ noetherian, we have that
every object $M\in {}_{H}{\mathcal {YD}}^H$ which is finitely generated as an
$H$-module is a direct sum in ${}_{H}{\mathcal {YD}}^H$ of a family of 
simple subobjects
that are finitely generated as $H$-modules.   
\end{corollary}

\begin{proof} 
Let $N$ be a subobject of $M$ in ${}_{H}{\mathcal {YD}}^H$. 
Then $M/N$ is finitely generated as an $H$-module and we have 
an exact sequence 
\begin{equation}\eqlabel{3.2.1}
0\rightarrow N\rightarrow M \rightarrow M/N \rightarrow 0
\end{equation}
in ${}_{H}{\mathcal {YD}}^H$.
$N$ is finitely generated as $H$-module, since 
$H$ is noetherian, so it follows from \prref{3.1} that $M/N$ and $N$ are projective in
${}_{H}{\mathcal {YD}}^H$, hence the sequence \equref{3.2.1} splits in
${}_{H}{\mathcal{YD}}^H$.
\end{proof}

Take $M\in {}_{H}{\mathcal {YD}}^H$ and $V$ a right $H$-subcomodule of $M$.
We will set
$$HV=\{\sum_{i\in I}a_iv_i~|~ a_i \in H, v_i \in V, \quad \hbox{where} \quad I \quad
\hbox{is a finite set}\}$$
$HV$ is a subobject of $M$ in ${}_{H}{\mathcal {YD}}^H$; the $H$-action and $H$-coaction
on
$HV$ are given by 
\begin{eqnarray*}
h(\sum_{i \in I}a_iv_i)&=& \sum_{i \in I}ha_iv_i\\
\rho(\sum_{i\in I}a_iv_i)&=&\sum_{i\in I}{(a_i)}_2{(v_i)}_0 \otimes {(a_i)}_3
{(v_i)}_1S^{-1}({(a_i)}_1)
\end{eqnarray*}

\begin{corollary}\colabel{3.3}
Let $H$ be commutative and noetherian. Assume that ${}_{H}{\mathcal {YD}}^H$
satisfies condition ($\dagger$), and that the functor $(-)^{{\rm co}H}$ from
$_{H}{\mathcal {YD}}^H$ to ${\mathcal M}$ is exact. Then $M\in
{}_{H}{\mathcal {YD}}^H$ is a direct sum in
${}_{H}{\mathcal {YD}}^H$ of a family of simple subobjects 
that are finitely generated as
$H$-modules. Therefore $M$ is a semisimple object in ${}_{H}{\mathcal {YD}}^H$ and
${}_{H}{\mathcal {YD}}^H$ is a semisimple category.
\end{corollary}

\begin{proof} 
Every $m\in M$ is contained in a finite-dimensional $H$-subcomodule $V_m$ of $M$,
see e.g. \cite[5.1.1]{DascalescuNR}. Then $HV_m$ is
finitely generated as $H$-module, and, by \coref{3.2},
each $HV_m$ is a direct sum of a family of simple subobjects of $HV_m$ 
(and of $M$) in ${}_{H}{\mathcal {YD}}^H$, which are finitely generated as an
$H$-module. Consequently each $m\in M$ is contained in a simple object which
is finitely generated as an $H$-module, so $M$ is a sum of simple objects
finitely generated as an $H$-module. The sum is a direct sum since the intersection
of two simple objects is trivial.
\end{proof}

\begin{corollary}\colabel{3.4}
Let $H$ be commutative, noetherian (in particular: finite dimensional), 
semisimple and cosemisimple. 
Then each $M\in {}_{H}{\mathcal {YD}}^H$ is a direct sum in ${}_{H}{\mathcal
{YD}}^H$ of a family of simple subobjects of $M$ finitely generated as $H$-modules.
Hence $M$ is semisimple in ${}_{H}{\mathcal {YD}}^H$ and ${}_{H}{\mathcal
{YD}}^H$ is a semisimple category.
\end{corollary}

\begin{proof} The cosemisimplicity of $H$ implies that the functor 
$$(-)^{{co}H}:\ {\mathcal M}^H\to {\mathcal M}$$
is exact, and, a fortiori
$$(-)^{{co}H}:\ {}_H{\mathcal YD}^H\to {\mathcal M}$$
is exact.
\end{proof}

Take $M,N\in {}_{H}{\mathcal L}^H$, with $M$ 
finitely generated as an $H$-module. By \prref{1.8} and \reref{1.9},
${}_H\Hom(M , N)\in {\mathcal M}^H$, and we can study the semisimplicity
of ${}_{H}{\mathcal L}^H$. We will say that ${}_{H}{\mathcal L}^H$
satisfies condition ($\dagger$) if the functor
$${}_H\Hom(M,-):\ {}_{H}{\mathcal L}^H\to {\mathcal M}^H$$
is exact for every $H$-finitely generated $M\in {}_{H}{\mathcal L}^H$.
The previous results of this Section then remain true after we replace
the category of Yetter-Drinfeld modules by Long dimodules, and without the
assumption that $H$ is commutative. We state the results without proof.

\begin{proposition}\prlabel{3.6}
Assume that ${}_{H}{\mathcal L}^H$ satisfies condition ($\dagger$) and that the
functor 
$$(-)^{{\rm co}H}:\ {\mathcal M}^H\to {\mathcal M}$$
is exact. Then every $H$-finitely generated $M\in {}_{H}{\mathcal L}^H$
is a projective object in ${}_{H}{\mathcal L}^H$.
\end{proposition}

\begin{corollary}\colabel{3.7}
Let $H$ be left noetherian, and assume that the conditions of \prref{3.6}
are satisfied. Then every $H$-finitely generated $M\in {}_{H}{\mathcal L}^H$
is a direct sum in ${}_{H}{\mathcal L}^H$ of a family of simple subobjects of $M$
that are finitely generated as $H$-modules. ${}_{H}{\mathcal L}^H$ is a semisimple
category.  
\end{corollary}

\begin{corollary}\colabel{3.9}
Let $H$ be left noetherian (in particular: finite dimensional), 
semisimple and cosemisimple. 
Then each $M\in {}_{H}{\mathcal L}^H$ is a direct sum in ${}_{H}{\mathcal L}^H$
of a family of simple subobjects of $M$ that are finitely generated as $H$-modules.
Hence
$M\in {}_{H}{\mathcal L}^H$ is semisimple and ${}_{H}{\mathcal L}^H$ is a
semisimple category.
\end{corollary}

\end{document}